\documentclass[reqno,12pt]{amsproc}
\usepackage{hyperref}
\usepackage{url}
\usepackage{amssymb}
\newtheorem{theo}{Theorem}
\newtheorem{rem}{Remark}

\newtheorem{probl}{Problem}

\newtheorem{prop}{Proposition}
\newtheorem{cor}{Corollary}
\newtheorem{df}{Definition}

\newcommand\eps\varepsilon
\newcommand\ph\varphi
\newcommand\kap\varkappa

\usepackage{enumerate}
\usepackage{graphicx}
\usepackage{float}
\usepackage{yfonts}
\usepackage{mathrsfs}
\usepackage{amsthm}
\usepackage{amsmath}
\usepackage{amscd}
\usepackage{xcolor}
\begin{document}\title
{Lecture Notes in Integral Invariants and Hamiltonian Systems}
\author[Oleg Zubelevich]{Oleg Zubelevich\\Faculty of Mechanics and Mathematics,\\
Lomonosov Moscow State University,\\ 1 Leninskie Gory, Moscow 119991, Russia}
\email{oezubel@gmail.com}\begin{abstract}
In this methodological review, we discuss the fundamental concepts of the theory of integral invariants. This theory originated with Poincare and Cartan \cite{Koz, Kart} and was further developed by  Kozlov \cite{int_K}. We demonstrate how the core ideas of this theory link diverse fields of mathematical physics, such as Hamiltonian dynamics, optics, and hydrodynamics. Particular attention is paid to results that are rarely expounded in standard textbooks.
\end{abstract}
\maketitle
\maketitle\tableofcontents
\section{Invariant Differential Forms}
Consider a smooth dynamical system
\begin{equation}\label{zsdgg}
\dot{x} = v(x),
\end{equation}
defined on the smooth manifold $M$ with local coordinates $x = (x^1, \dots, x^m)^T,\quad v(x)\in T_xM$.
\begin{rem}\label{xdfgg5ttyy8ppoolgdfgh}Throughout this text, all objects are assumed to be sufficiently smooth for the respective expressions to be well-defined, and all applied theorems are assumed to have their conditions satisfied. All integrals are assumed to exist. Whenever Stokes’ theorem is applied, we assume that the necessary regularity and boundary conditions are met. Etc.\end{rem}

We denote the phase flow of this system by $g^t: M \to M$, which satisfies:
$$\frac{d}{dt}g^t({x}) = v(g^t({x})), \quad g^0({x}) = {x}.$$

Recall the group property $g^{t+s} = g^t \circ g^s$.

\begin{df}[\cite{nov}]\label{sdg000d}
Consider a differential $k$-form
$$
\omega = \sum_{i_1 < \dots < i_k} \omega_{i_1 \dots i_k}(x) dx^{i_1} \wedge \dots \wedge dx^{i_k}.
$$
The Lie derivative of $\omega$ along the vector field $v$, denoted by $L_v \omega$, is defined as the following differential operator:
$$
L_v \omega = \left. \frac{d}{dt} \right|_{t=0} g^t_* \omega = \left. \frac{d}{dt} \right|_{t=0} \sum_{i_1 < \dots < i_k} \omega_{i_1 \dots i_k}(g^t(x)) d(g^t(x))^{i_1} \wedge \dots \wedge d(g^t(x))^{i_k}.$$
\end{df}
The operator $i_v \omega$ is called the \textit{interior product} (or \textit{contraction}) of the form $\omega$ with the vector field $v$. It maps a $k$-form ($k \in \mathbb{N}$) to a $(k-1)$-form; at each point $x \in M$, it is defined by
$$(i_v\omega)(u_1, \ldots, u_{k-1}) = \omega(v(x), u_1, \ldots, u_{k-1}),$$
where $u_1, \ldots, u_{k-1} \in T_xM$.

This operator satisfies the following antiderivation property:
\begin{equation}\label{sswe4rfv}
i_v(\omega \wedge \nu) = (i_v \omega) \wedge \nu + (-1)^k \omega \wedge (i_v \nu),
\end{equation}
where $k$ is the degree of $\omega$.

For a basis $k$-form, the following identity holds:
$$
i_v(dx^{i_1} \wedge \ldots \wedge dx^{i_k}) = \sum_{s=1}^k (-1)^{s+1} v^{i_s} dx^{i_1} \wedge \ldots \wedge \widehat{dx^{i_s}} \wedge \ldots \wedge dx^{i_k},
$$
where the hat symbol $\widehat{\ldots}$ denotes that the corresponding factor is omitted from the wedge product.

Recall the Cartan  magic formula:
\begin{equation}\label{zsdg00v}L_v \omega = d i_v \omega + i_v d \omega.\end{equation}

Recall that the exterior derivative satisfies the same antiderivation property as $i_v$:
$$
d(\omega \wedge \nu) = (d\omega) \wedge \nu + (-1)^k \omega \wedge (d\nu).
$$

Furthermore, the following formula holds:
\begin{equation}\label{Psfgb66} L_u i_v - i_v L_u = i_{[u, v]}. \end{equation}

\begin{probl}
Show that
\begin{align}
L_v\big(\omega_idx^i\big)&=\Big(\frac{\partial\omega_i}{\partial x^s}v^s+\frac{\partial v^s}{\partial x^i}\omega_s\Big)dx^i,\label{sdgftth}\\
L_v\Big(\rho(x)dx^1\wedge\ldots\wedge dx^m\Big)&=\frac{\partial (\rho v^i)}{\partial x^i}dx^1\wedge\ldots\wedge dx^m,\label{zsdf431qqw}\\
L_v f&=v^i\frac{\partial f}{\partial x^i},\quad f:M\to \mathbb{R},\nonumber\\
L_v(\omega\wedge\nu)&=(L_v\omega)\wedge\nu+\omega\wedge(L_v\nu),\nonumber\\
L_vd&=dL_v.\label{sdfgtymmm}
\end{align}
\end{probl}
\begin{probl} Using $(\ref{zsdg00v})$ and $(\ref{Psfgb66})$, verify the following identity:
$$L_{[u,v]}=L_uL_v-L_vL_u.$$\end{probl}

\begin{df}A $k$-form $\omega$ is called an integral invariant of system \eqref{zsdgg} if and only if $L_v \omega = 0$.

A $k$-form $\omega$ is called a relative integral invariant of system \eqref{zsdgg} if there exists a $(k-1)$-form $\Omega$ such that $L_v \omega = d\Omega$.
\end{df}
\begin{rem}\label{xfgb660} From formula \eqref{sdfgtymmm}, it follows that if the form $\omega$ is a relative integral invariant, then the form $d\omega$ is an (absolute) integral invariant.\end{rem}

Let $\Sigma \subset M$ denote a smooth $k$-dimensional submanifold of $M$.
\begin{theo}\label{dsfsfd}The following identity holds:
\begin{equation}\label{vbvgbfg}
\left. \frac{d}{dt} \right|_{t=0} \int_{g^t(\Sigma)} \omega = \int_\Sigma L_v \omega.
\end{equation}
\end{theo}
\subsubsection*{Proof} We prove this identity for the case where $\Sigma$ is covered by a single coordinate chart.

Let $u: D \to M$ denote the corresponding embedding, where $D \subset \mathbb{R}^k$ is an open domain such that $u(D) = \Sigma$.

Thus, by applying the change of variables formula, we obtain:
$$\int_{g^t(\Sigma)}\omega=\int_D(g^t\circ u)_*\omega=\int_D(u_*\circ (g^t)_*)\omega,$$ and
$$\frac{d}{dt}\Big|_{t=0}\int_D(u_*\circ (g^t)_*)\omega=\int_Du_*\frac{d}{dt}\Big|_{t=0}(g^t)_*\omega=\int_Du_*(L_v\omega).$$

The theorem is proved.

\begin{cor}Formula (\ref{vbvgbfg}) implies
\begin{equation}\label{sdfgbvf}\frac{d}{dt}\int_{g^t(\Sigma)}\omega=
\frac{d}{ds}\Big|_{s=0}\int_{g^{t+s}(\Sigma)}\omega=\frac{d}{ds}\Big|_{s=0}\int_{g^s\big(g^t(\Sigma)\big)}\omega=
\int_{g^t(\Sigma)}L_v\omega.\end{equation}\end{cor}

\begin{theo}\label{sdfr5678} If $\omega$ is an integral invariant, then for any admissible $k$-dimensional submanifold $\Sigma \subset M$ (see Remark \ref{xdfgg5ttyy8ppoolgdfgh}), the integral
$$
\int_{g^t(\Sigma)} \omega
$$
is independent of $t$. The converse is also true.\end{theo}
Indeed, this follows directly from \eqref{sdfgbvf}.

\begin{theo}\label{sdgffrt} Assume that $\omega$ is a relative integral invariant. Then for any compact $k-$dimensional submanifold $\Sigma\subset M,\quad \partial \Sigma=\emptyset$ the following integral $$\int_{g^t(\Sigma)}\omega$$ does not depend on $t$. \end{theo}
Indeed, from formula \eqref{sdfgbvf} and Stokes' theorem, we have:
$$
\frac{d}{dt} \int_{g^t(\Sigma)} \omega = \int_{g^t(\Sigma)} L_v \omega = \int_{g^t(\Sigma)} d\Omega = 0.
$$
This holds because $\partial g^t(\Sigma) = \emptyset$.

\begin{probl}Assume that the system \eqref{zsdgg} has two integral invariants$$\omega_i = \rho_i(x) dx^1 \wedge \ldots \wedge dx^m, \quad i=1, 2, \quad \rho_2 \neq 0.$$Show that the ratio $\rho_1 / \rho_2$ is a first integral of the system \eqref{zsdgg}.\end{probl}

\begin{probl}
Show that if $L_v\omega = 0$, then the zero set of the form,
$$Q = \{x \in M \mid \omega(x) = 0\},$$
is invariant under the corresponding phase flow $g^t$, i.e., $g^t(Q) = Q$ for all $t$.
\end{probl}

\section{Invariant Differential Forms of Systems with First Integrals}
\subsection{The Case $\dim M=2$}\label{dfg56}
If $m=2$ and $q$ is an invariant $2$-form of system (\ref{zsdgg}):
$$L_v q = 0,$$
then $i_v q$ is a closed form. Indeed, by Cartan's formula:$$L_v q = di_v q + i_v dq.$$
Having $dq=0$, we obtain $di_v q = 0$.

Thus, locally there exists a function $f$ such that:$$df = i_v q.$$
It follows that:$$L_v f = i_v df = i_v i_v q = 0.$$
Consequently, the function $f$ is a first integral.

 \subsection{Invariant Forms on Level Sets of First Integrals}

Assume that system (\ref{zsdgg}) has a first integral $F$ such that:
$$L_v F = 0, \quad dF \neq 0,$$
and an invariant $m$-form $\omega:\quad L_v \omega = 0$.

Without loss of generality, we can locally assume that:
$$\frac{\partial F}{\partial x^m} \neq 0.$$
Then $\omega$ can be written as:
$$\omega = \rho(x) dx^1 \wedge \dots \wedge dx^m = \rho(x) \left( \frac{\partial F}{\partial x^m} \right)^{-1} dx^1 \wedge \dots \wedge dx^{m-1} \wedge dF.$$
In invariant terms, this means that the form $\omega$ is represented as the following exterior product:$$\omega = \lambda \wedge dF,$$
where $\lambda$ is an $(m-1)$-form defined up to an additive form $\gamma$ such that $\gamma \wedge dF = 0$.

\begin{prop}The restricted differential form$$\lambda|_{Z}, \quad Z = \{F = \text{const}\},$$is an integral invariant for the restriction of system (\ref{zsdgg}) to the manifold $Z$.\end{prop}
Indeed,
$$L_v \omega = (L_v \lambda) \wedge dF
+ \lambda \wedge (L_v dF), \quad L_v(dF) = dL_v F = 0.$$
Therefore, we have $(L_v \lambda) \wedge dF = 0$.

Let $e_1, \dots, e_m$ be a basis in $T_x M$ such that $e_1, \dots, e_{m-1}$ is a basis in $T_x Z$.

We consequently obtain:$$\big((L_v \lambda) \wedge dF\big)(e_1, \dots, e_m) = c (L_v \lambda)(e_1, \dots, e_{m-1}) \cdot dF(e_m) = 0,$$
where $c$ is a non-zero constant related to the definition of the exterior product.

Since $dF(e_m) \neq 0$, it follows that $(L_v \lambda)(e_1, \dots, e_{m-1}) = 0$.

The proposition is proved.

Taking into account the results of Section \ref{dfg56}, we see that if system (\ref{zsdgg}) has $m-2$ independent first integrals and an invariant $m$-form:
$$L_v \omega = 0, \quad \omega = \rho(x) dx^1 \wedge \dots \wedge dx^m, \quad \rho(x) > 0,$$
then this system is integrable in closed form.

 \section{The Lie Derivative in the Nonautonomous case}
Consider a vector field$$v=(v^1,\ldots, v^m)(t,x), \quad (t,x) \in \tilde M = (t_1,t_2) \times M.$$
We use $G^t_{t_0}: M \to M$ to denote the shift along the trajectories of the system
\begin{equation}\label{zsdfff}\dot x = v(t,x), \quad \frac{dG^t_{t_0}(x)}{dt} = v(t,G^t_{t_0}(x)), \quad G^{t_0}_{t_0}(x) = x, \quad G_{t_1}^t \circ G_{t_0}^{t_1} = G_{t_0}^t.\end{equation}

Let $g^\tau: \tilde M \to \tilde M$ be the flow of the system
$$\frac{dz}{d\tau} = \tilde v(z), \quad z = (t, x^1, \ldots, x^m)^T, \quad \tilde v = (1, v^1, \ldots, v^m)^T.$$
Here, $\tilde M$ denotes the extended phase space of (\ref{zsdfff}).

The following equality holds:
$$g^\tau(t_0,x)=(t_0+\tau,G_{t_0}^{t_0+\tau}(x)).$$

Consider a form
$$\omega = \sum_{i_1 < \ldots < i_k} \omega_{i_1 \ldots i_k}(t,x) dx^{i_1} \wedge \ldots \wedge dx^{i_k}$$
and fix local coordinates.
We can consider $\omega$ as a form on the manifold $\tilde M$ or as a form on the manifold $M$. In the latter case, $t$ is regarded as a parameter.

Introduce the following notations:
$$
\frac{\partial\omega}{\partial t}=\sum_{i_1<\ldots<i_k}\frac{\partial\omega_{i_1\ldots i_k}(t,x)}{\partial t} dx^{i_1}\wedge\ldots\wedge dx^{i_k},\quad
\frac{\partial v}{\partial t}=\left(\frac{\partial v^1}{\partial t},\ldots,\frac{\partial v^m}{\partial t}\right);
$$
and
\begin{align}
d &= d_t + d_x, \quad d_t\omega = \sum_{i_1<\ldots<i_k}\frac{\partial\omega_{i_1\ldots i_k}(t,x)}{\partial t}dt\wedge dx^{i_1}\wedge\ldots\wedge dx^{i_k}; \nonumber \\
d_x\omega &= \sum_{i_1<\ldots<i_k}\frac{\partial\omega_{i_1\ldots i_k}}{\partial x^l}dx^l\wedge dx^{i_1}\wedge\ldots\wedge dx^{i_k}, \nonumber \\
L_v\omega &:= d_x i_v \omega + i_v d_x \omega. \nonumber
\end{align}
For instance, we have $d_t\omega = dt \wedge \frac{\partial\omega}{\partial t}$.
\begin{theo}\label{xfbbb}
The following formula holds:
$$L_{\tilde v}\omega = \frac{\partial\omega}{\partial t} + L_v\omega + (dt) \wedge i_{\frac{\partial v}{\partial t}}\omega.$$\end{theo}
In particular, one finds that
$$L_{\tilde v}\omega\Big|_{t=\mathrm{const}} = \frac{\partial\omega}{\partial t} + L_v\omega.$$
\subsubsection*{Proof of Theorem  \ref{xfbbb}}Let us split the vector field into the following summands:
$$ \tilde v = e + v_*, \quad e = (1, 0, \ldots, 0), \quad v_* = (0, v^1, \ldots, v^m). $$
Correspondingly, $L_{\tilde v} = L_e + L_{v_*}.$

Observe that
$$ L_e\omega = \frac{\partial\omega}{\partial t}, \quad i_{v_*}\omega = i_v\omega. $$
By the Cartan magic formula, we get:
\begin{equation}\label{xcvtt6}
L_{v_*}\omega = i_{v_*}d_x\omega + i_{v_*}d_t\omega + d_xi_{v_*}\omega + d_ti_{v_*}\omega.
\end{equation}
Since
$$ i_{v_*}d_x\omega = i_{v}d_x\omega, \quad d_xi_{v_*}\omega = d_xi_{v}\omega, $$
equality (\ref{xcvtt6}) takes the form:
$$ L_{v_*}\omega = L_v\omega + i_{v_*}d_t\omega + d_ti_{v_*}\omega. $$
By formula (\ref{sswe4rfv}), it follows that:
$$ i_{v_*}d_t\omega = i_{v_*}\left((dt) \wedge \frac{\partial\omega}{\partial t}\right) = -(dt) \wedge i_v\frac{\partial\omega}{\partial t}, $$
and
$$ d_ti_{v_*}\omega = (dt) \wedge \frac{\partial i_v\omega}{\partial t} = (dt) \wedge \left(i_{\frac{\partial v}{\partial t}}\omega + i_v\frac{\partial\omega}{\partial t}\right). $$
This proves the theorem.

\begin{theo}\label{sxdgff5jtklgt000}For any admissible $k$-dimensional submanifold $A \subset M$, one has:
\begin{equation}\label{xfgg597uuji}
\frac{d}{dt} \int_{G^t_{t_0}(A)} \omega(t,\cdot) = \int_{G^t_{t_0}(A)} \left( \frac{\partial\omega}{\partial t}(t,\cdot) + L_{v(t,\cdot)}\omega(t,\cdot) \right).
\end{equation}\end{theo}
\subsubsection*{Proof}Introduce the manifold $A_{t} = \{t\} \times A \subset \tilde M$, where $g^\tau(A_{t}) = \{t+\tau\} \times G^{t+\tau}_{t}(A) \subset \tilde M$.

Using Theorem \ref{xfbbb} and formula \ref{sdfgbvf}, we calculate:
\begin{align}
\frac{d}{d\tau} \int_{g^\tau(A_{t_0})} \omega &= \int_{g^\tau(A_{t_0})} L_{\tilde v} \omega = \int_{G^{t_0+\tau}_{t_0}(A)} (L_{\tilde v} \omega) \big|_{t=t_0+\tau} \nonumber \\
&= \int_{G^{t_0+\tau}_{t_0}(A)} \left( \frac{\partial\omega}{\partial t} + L_v \omega \right) \bigg|_{t=t_0+\tau}. \nonumber
\end{align}

To finish the proof, it remains to observe that
$$ \int_{g^\tau(A_{t_0})} \omega = \int_{G^{t_0+\tau}_{t_0}(A)} \omega \big|_{t=t_0+\tau}. $$

The theorem is proved.

\begin{rem}After the change of variables $x = G^t_{t_0}(y)$, formula (\ref{xfgg597uuji}) takes the form:
$$
\frac{d}{dt} \int_{A} (G^t_{t_0})_* \omega(t,\cdot) =
\int_{A} \frac{d}{dt} (G^t_{t_0})_* \omega(t,\cdot)
= \int_{A} (G^t_{t_0})_* \left( \frac{\partial\omega}{\partial t}(t,\cdot) + L_{v(t,\cdot)}\omega(t,\cdot) \right).
$$
Since $A$ is an arbitrary $k-$dimensional manifold, we have:
\begin{equation}\label{zdggty}
\frac{d}{dt} (G^t_{t_0})_* \omega(t,\cdot) = (G^t_{t_0})_* \left( \frac{\partial\omega}{\partial t}(t,\cdot) + L_{v(t,\cdot)}\omega(t,\cdot) \right).
\end{equation}\end{rem}

\begin{theo}\label{zxcvdv001}
1) Let $A \subset M$ be an admissible $k$-dimensional submanifold and assume that
$$ \frac{\partial\omega}{\partial t} + L_v\omega = 0. $$
Then we have
$$ \int_{G^t_{t_0}(A)} \omega(t,\cdot) = \int_A \omega(t_0,\cdot). $$

2) Suppose that there exists a $(k-1)$-form
$$ \Omega = \sum_{i_1 < \ldots < i_{k-1}} \Omega_{i_1 \ldots i_{k-1}}(t,x) dx^{i_1} \wedge \ldots \wedge dx^{i_{k-1}} $$
such that
$$ \frac{\partial\omega}{\partial t} + L_v\omega = d_x \Omega, $$
and $A \subset M$ is a compact submanifold without boundary ($\partial A = \emptyset$). Then we have
$$ \int_{G^t_{t_0}(A)} \omega(t,\cdot) = \int_A \omega(t_0,\cdot). $$
\end{theo}

Theorem \ref{zxcvdv001} follows from Theorem \ref{sxdgff5jtklgt000} in the same manner as Theorems \ref{sdgffrt} and \ref{sdfr5678} follow from formula (\ref{sdfgbvf}).

The theorem is proved.

\section{A Remark on Solutions to a PDE}Formula (\ref{zdggty}) provides a method for solving the following Cauchy problem:
\begin{equation}\label{xvb77uuii90}
\frac{\partial\omega}{\partial t}(t,x) + L_{v(t,x)}\omega(t,x) = 0, \quad \omega\mid_{t=t_0} = \hat\omega(x).
\end{equation}
Indeed, from this formula it follows that
$$\frac{d}{dt} (G^t_{t_0})_* \omega(t,\cdot) = 0.$$
Since $(G^t_{t_0})_* \omega(t,\cdot)$ does not depend on $t$, we obtain
$$(G^t_{t_0})_* \omega(t,\cdot) = \hat\omega, \quad \text{whence} \quad \omega(t,x) = \big((G^t_{t_0})^{-1}\big)_* \hat\omega.$$

As an example, consider a linear vector field $v(t,x)=A(t)x,\quad x\in\mathbb{R}^m$, where $A$ is an $m \times m$ matrix. The shift $G_0^t$ is given by the formula
$$G_0^t(x)=X(t)x,$$
where $X$ is the fundamental matrix satisfying
$$\dot X=A(t)X, \quad X(0)=I.$$

Consider a differential form
$$\hat\omega = \hat\rho dx^1 \wedge \ldots \wedge dx^m,$$
where $\hat\rho \neq 0$ is a constant. Then the form
$$\omega(t) = \rho(t) dx^1 \wedge \ldots \wedge dx^m = \big((G^t_{0})^{-1}\big)_* \hat\omega = \det X^{-1}(t) \hat\rho dx^1 \wedge \ldots \wedge dx^m$$
is a solution to the Cauchy problem (\ref{xvb77uuii90}) (with $t_0=0$). By formula (\ref{zsdf431qqw}), this is equivalent to saying that
$$\rho_t + \mathrm{div}(\rho v) = 0, \quad \rho|_{t=0} = \hat\rho.$$
Observe also that since $\rho(t)$ does not depend on $x$, we have
$$\mathrm{div}(\rho v) = \rho(t) \mathrm{tr} A(t).$$
It follows that
$$\rho(t) = \hat\rho e^{-\int_0^t \mathrm{tr} A(s) ds} = \det X^{-1}(t) \hat\rho.$$
We arrive at the celebrated Liouville's formula:
$$\det X(t) = e^{\int_0^t \mathrm{tr} A(s) ds}.$$
\section{Applications to Hydrodynamics}
In this section, we set $M = \mathbb{R}^3$ and let $(x^1, x^2, x^3)$ be the standard right-handed Euclidean frame. We also adopt the convention $d = d_x$.

Let $\mathbf{A}(t,x) = (A_i)$ and $\mathbf{B}(t,x) = (B_i)$ be vector fields in $\mathbb{R}^3$, and let $f(t,x)$ be a scalar function. There are standard correspondences:
$$ f \mapsto \omega^3_f = f \, dx^1 \wedge dx^2 \wedge dx^3 $$
and
$$ \mathbf{A} \mapsto \omega^1_{\mathbf{A}} = A_i dx^i, \quad \mathbf{A} \mapsto \omega^2_{\mathbf{A}} = A_1 dx^2 \wedge dx^3 + A_2 dx^3 \wedge dx^1 + A_3 dx^1 \wedge dx^2. $$
These correspondences can be expressed in terms of the Hodge star operator.

The following formulas are verified by direct calculation:
$$ df = \omega^1_{\mathrm{grad}\,f}, \quad d\omega^1_{\mathbf{A}} = \omega^2_{\mathrm{curl}\,\mathbf{A}}, \quad d\omega^2_{\mathbf{A}} = \omega^3_{\mathrm{div}\,\mathbf{A}} $$
and
$$ i_{\mathbf{B}}\omega^1_{\mathbf{A}} = (\mathbf{A}, \mathbf{B}), \quad i_{\mathbf{B}}\omega^2_{\mathbf{A}} = \omega^1_{\mathbf{A} \times \mathbf{B}}, \quad i_{\mathbf{B}}\omega^3_{f} = f \omega^2_{\mathbf{B}}. $$

These identities imply the following theorem.
\begin{theo}\label{xcxcxc}
The following formulas hold:
\begin{align}
\frac{\partial \omega^1_{\mathbf{A}}}{\partial t} + L_v \omega^1_{\mathbf{A}} &=
\omega^1_{\frac{\partial \mathbf{A}}{\partial t} + (\mathrm{curl}\,\mathbf{A}) \times v} + d(v, \mathbf{A}), \nonumber \\
\frac{\partial \omega^2_{\mathbf{A}}}{\partial t} + L_v \omega^2_{\mathbf{A}} &=
\omega^2_{\frac{\partial \mathbf{A}}{\partial t} + \mathrm{curl}\,(\mathbf{A} \times v) + v\,\mathrm{div}\,\mathbf{A}}, \nonumber \\
\frac{\partial \omega^3_{f}}{\partial t} + L_v \omega^3_{f} &=
\omega^3_{\frac{\partial f}{\partial t} + \mathrm{div}\,(f v)}. \nonumber
\end{align}
\end{theo}

\begin{theo}\label{sxdfgggcde44}
Let $\gamma$, $\Sigma$, and $D \subset \mathbb{R}^3$ be a closed curve, a two-dimensional surface, and a domain, respectively.

1) If
\begin{equation}\label{zaaadvff}
\frac{\partial \mathbf{A}}{\partial t} + (\mathrm{curl}\,\mathbf{A}) \times v = \mathrm{grad}\,\psi, \quad \psi = \psi(t,x),
\end{equation}
then
$$ \int_{G_{t_0}^t(\gamma)} \omega^1_{\mathbf{A}(t,\cdot)} = \int_{\gamma} \omega^1_{\mathbf{A}(t_0,\cdot)}; $$

2) If
\begin{equation}\label{zxcvfgtt6}
\frac{\partial \mathbf{A}}{\partial t} + \mathrm{curl}\,(\mathbf{A} \times v) + v\,\mathrm{div}\,\mathbf{A} = 0,
\end{equation}
then
$$ \int_{G_{t_0}^t(\Sigma)} \omega^2_{\mathbf{A}(t,\cdot)} = \int_{\Sigma} \omega^2_{\mathbf{A}(t_0,\cdot)}; $$

3) If
\begin{equation}\label{axdervvv}
\frac{\partial f}{\partial t} + \mathrm{div}\,(f v) = 0,
\end{equation}
then
$$ \int_{G_{t_0}^t(D)} \omega^3_{f(t,\cdot)} = \int_{D} \omega^3_{f(t_0,\cdot)}. $$
\end{theo}
This theorem is a consequence of Theorems \ref{zxcvdv001} and \ref{xcxcxc}.

\begin{probl}Prove the equality
$$ \mathrm{curl}\,(\mathbf{A} \times \mathbf{B}) = -[\mathbf{A}, \mathbf{B}] + \mathbf{A} \, \mathrm{div}\,\mathbf{B} - \mathbf{B} \, \mathrm{div}\,\mathbf{A}. $$
We use square brackets to denote the commutator of the vector fields:
$$ [\mathbf{A}, \mathbf{B}] = -\frac{\partial \mathbf{A}}{\partial x}\mathbf{B} + \frac{\partial \mathbf{B}}{\partial x}\mathbf{A}. $$
Hint: apply  formula (\ref{Psfgb66})   to the form $\omega=dx^1\wedge dx^2\wedge dx^3.$
\end{probl}

Let $v(t,x)$ be a velocity field of an ideal fluid subject to potential forces, and assume its density $\rho$ depends on pressure $p$ only (a barotropic fluid).

Under these assumptions, the equation of fluid's motion has the form:
\begin{equation}\label{xfb--0}
\frac{\partial v }{\partial t}+\boldsymbol \omega\times  v=
-\mathrm{grad}\,\mathcal P-\mathrm{grad}\,V-\mathrm{grad}\,\frac{|v|^2}{2},\quad \boldsymbol \omega=\mathrm{curl}\, v,
\end{equation}
where $-\mathrm{grad}\,V=\mathbf F$ is a force per unit mass and $\mathrm{grad}\,\mathcal P=\frac{1}{\rho}\mathrm{grad}\,p$.

So that equation (\ref{zaaadvff}) holds with $\mathbf{A} = v$.

Taking the $\mathrm{curl}$ of both sides of (\ref{xfb--0}) we obtain:
$$\frac{\partial \boldsymbol \omega}{\partial t}+\mathrm{curl}\,(\boldsymbol \omega\times  v)=0.$$

Consequently the equation (\ref{zxcvfgtt6}) holds with $\mathbf{A} = \boldsymbol \omega$.

The famous Helmholtz and Kelvin theorems follow from Theorem \ref{sxdfgggcde44}. Furthermore, for $f=\rho$ the equation (\ref{axdervvv}) corresponds to the continuity equation.

For further reading see \cite{coch}.

\section{ Darboux Theorem \cite{Darboux}, \cite{Darboux1}}\label{zdsf4rt}
In this section, we apply the developed theory to a more complex case.

Let $D \subset \mathbb{R}^{2m} = \{x = (x^1, \dots, x^{2m})\}$ be an open neighbourhood of the origin.

Suppose a differential form
$$ \omega = \sum_{i<j} \omega_{ij}(x) dx^i \wedge dx^j $$
defined in $D$ is closed and non-degenerate:
$$ d\omega = 0, \quad \det(\omega_{ij}(x)) \ne 0, \quad x \in D. $$
\begin{theo}\label{zadgv5ty78}
In some open neighbourhood $D' \subset D$ of the origin, there exist local coordinates $y = (y^1, \dots, y^{2m})$ such that the form $\omega$ has constant coefficients:
$$ \omega = \sum_{i<j} \omega'_{ij} dy^i \wedge dy^j, \quad \omega'_{ij} = \text{const}_{ij}. $$
\end{theo}
From linear algebra, it is known that there exists a linear transformation $y \mapsto (q^1, \dots, q^m, p_1, \dots, p_m)$ such that
$$ \omega =  dp_i \wedge dq^i. $$
\subsubsection*{Proof of Theorem \ref{zadgv5ty78}}
Introduce a constant form
$$ \omega_1 = \omega|_{x=0} = \sum_{i<j} \omega_{ij}(0) dx^i \wedge dx^j, \quad d\omega_1 = 0. $$
Construct a family of differential forms as follows:
$$ \Omega(t,x) = t\omega_1 + (1-t)\omega = \sum_{i<j} \Omega_{ij}(t,x) dx^i \wedge dx^j, \quad t \in [0,1], \quad d_x \Omega = 0. $$
These forms are well-defined in $D$.

Each form $\Omega(t,0) = \omega|_{x=0}$ is non-degenerate, so all the forms $\{\Omega(t,x)\}$ are non-degenerate in some neighbourhood of the origin. Since $d(\omega_1 - \omega) = 0$, by the Poincare lemma, there exists a $1$-form $\alpha$ such that
$$ \omega_1 - \omega = d\alpha, \quad \alpha =  \alpha_k(x) dx^k. $$
Let us choose $\alpha$ such that $\alpha|_{x=0} = 0$.

Define a time-dependent vector field $v(t,x)$ as follows:
$$ i_{v(t,x)} \Omega(t,x) = -\alpha, \quad v^i(t,x) = -\Omega^{ij}(t,x) \alpha_j(x). $$
Since $v(t,0) = 0$, the corresponding family of shifts $x \mapsto G^t_0(x)$ (along trajectories of the system $\dot x=v(t,x)$) is defined for small $|x|$ and for $t \in [0,1]$.

Using  formula (\ref{zdggty}), we obtain
$$ \frac{d}{dt} (G^t_0)_* \Omega(t,\cdot) = (G^t_0)_* \left( \frac{\partial\Omega}{\partial t}(t,\cdot) + L_{v(t,\cdot)} \Omega(t,\cdot) \right). $$
From Cartan's formula, we get
$$ L_{v(t,\cdot)} \Omega(t,\cdot) = i_{v(t,x)} d_x \Omega(t,x) + d_x i_{v(t,x)} \Omega(t,x) = -d\alpha. $$
On the other hand,
$$ \frac{\partial\Omega}{\partial t} = \omega_1 - \omega = d\alpha. $$
Thus, we have
$$ \frac{d}{dt} (G^t_0)_* \Omega(t,\cdot) = 0. $$
Consequently, $(G^1_0)_* \Omega(1,\cdot) = (G^0_0)_* \Omega(0,\cdot) = \omega$, and the mapping $G^1_0$ takes the form $\omega_1$ to the form $\omega$.

The theorem is proved.

\section{Integral Invariants of the Hamilton Equations}
Recall that a manifold $N$ with $\dim N = 2m$ is called a \textbf{symplectic manifold} if it is endowed with a non-degenerate closed 2-form, denoted by $\beta$. Non-degeneracy means that the equality $i_a\beta = 0$ implies $a = 0$.

Darboux's theorem (see Section \ref{zdsf4rt}) states that any point of $N$ belongs to a chart with local coordinates
$$z=(x^1,\ldots,x^m,p_1,\ldots, p_m)$$
such that
$$\beta=dp_i\wedge dx^i.$$
These coordinates are called symplectic or canonical.

Consider the Hamiltonian system with the Hamiltonian function $H=H(t,x,p)$:
\begin{equation}\label{sdfggg}
\dot x^i=\frac{\partial H}{\partial p_i},\quad \dot p_i=-\frac{\partial H}{\partial x^i},\quad t\in I=(t_1,t_2).
\end{equation}

Let $w$ stand for the corresponding vector field:
\begin{equation}\label{xdfsg500i}
w(t,z)=\Big(\frac{\partial H}{\partial p_1},\ldots,\frac{\partial H}{\partial p_m}, -\frac{\partial H}{\partial x^1},\ldots,-\frac{\partial H}{\partial x^m}\Big).
\end{equation}

Note that $w$ is uniquely determined by the equation
\begin{equation}\label{sdgf5txsw3}
i_w\beta=-d_zH.
\end{equation}
Therefore, this equation can be taken as an invariant definition of a Hamiltonian vector field.

Alongside system (\ref{sdfggg}), it is convenient to study its autonomous version in the extended phase space:
\begin{equation}\label{sxdfgyui}
\frac{dt}{d\tau}=1,\quad \frac{dx^i}{d\tau}=\frac{\partial H}{\partial p_i},\quad \frac{dp_i}{d\tau}=-\frac{\partial H}{\partial x^i}.
\end{equation}
The corresponding vector field is
$$\tilde w(t,z)=\Big(1, \frac{\partial H}{\partial p_1},\ldots,\frac{\partial H}{\partial p_m}, -\frac{\partial H}{\partial x^1},\ldots,-\frac{\partial H}{\partial x^m}\Big).$$
We use $\tilde N$ to denote the phase space of system (\ref{sxdfgyui}). The manifold $\tilde N$ is the extended phase space of (\ref{sdfggg}):
$$(t,z)\in\tilde N=I\times N.$$
Let $g^\tau:\tilde N\to \tilde N$ stand for the flow of (\ref{sxdfgyui}).

Let $G_{t_0}^t:N\to N$ denote the shift along the trajectories of (\ref{sdfggg}):
$$\frac{dG^t_{t_0}(z)}{dt}=w(t,G^t_{t_0}(z)),\quad G^{t_0}_{t_0}(z)=z.$$
Recall that
$$g^\tau\big((t_0,z)\big)=\big(t_0+\tau,G^{t_0+\tau}_{t_0}(z)\big).$$
Introduce a differential form on $\tilde N$:
$$\alpha=p_idx^i-Hdt.$$
This form is referred to as the Poincare--Cartan relative integral invariant. We justify this name below.
\begin{theo}\label{cccfg}
The following formula holds:
$$
i_{\tilde w}d\alpha=0.
$$
Conversely, if a vector field $u(t,z)$ satisfies the equality $i_u d\alpha=0$, then $u=\lambda(t,z)\tilde w$ for some scalar function $\lambda(t,z)$.
\end{theo}
\subsubsection*{Proof}The 2-form $d\alpha$ can be represented as follows:
$$
d\alpha = dp_i \wedge dx^i - dH \wedge dt = \pi_i \wedge \varkappa^i,
$$
where
$$ \pi_i = dp_i + \frac{\partial H}{\partial x^i} dt, \quad \varkappa^i = dx^i - \frac{\partial H}{\partial p_i} dt. $$
The value of the form $d\alpha$ on the vectors $\xi, \eta$ is calculated as
$$ d\alpha(\xi, \eta) = \pi_i(\xi) \varkappa^i(\eta) - \pi_i(\eta) \varkappa^i(\xi). $$
Since $\varkappa^i(\tilde w) = 0$ and $\pi_i(\tilde w) = 0$ by construction, we have
$$ i_{\tilde w} d\alpha = \pi_i(\tilde w) \varkappa^i(\cdot) - \pi_i(\cdot) \varkappa^i(\tilde w) = 0. $$
This proves the first part of the theorem.

To prove the second part, assume that $i_u d\alpha = 0$, which implies
$$ \pi_i(u) \varkappa^i(\cdot) - \pi_i(\cdot) \varkappa^i(u) = 0. $$
Note that the 1-forms $\pi_i, \varkappa^i$ for $i = 1, \ldots, m$ are linearly independent; therefore,
$$ \varkappa^i(u) = 0 \quad \text{and} \quad \pi_i(u) = 0. $$
Setting $u = (u_t, u_x, u_p)$, we obtain the relations
$$ u_{p_i} + \frac{\partial H}{\partial x^i} u_t = 0, \quad u_{x^i} - \frac{\partial H}{\partial p_i} u_t = 0. $$
Thus, $u = u_t \tilde w$, which completes the proof.

\begin{theo}\label{sdfgtaa}
The form $\alpha$ is a relative integral invariant of the system (\ref{sxdfgyui}):
$$
L_{\tilde w}\alpha = d\mathscr{F}, \quad \text{where} \quad \mathscr{F} = p_i \frac{\partial H}{\partial p_i} - H.
$$\end{theo}
\subsubsection*{Proof}This can be verified by direct calculation using formula (\ref{sdgftth}). However, it is more convenient to use the Cartan formula and Theorem \ref{cccfg}:
$$
L_{\tilde w}\alpha = d(i_{\tilde w}\alpha) + i_{\tilde w}d\alpha, \quad \text{where} \quad d(i_{\tilde w}\alpha) = d \mathscr{F}.
$$
The theorem is proved.

Note that if the Hamiltonian $H$ admits the Legendre transform with respect to the momenta $p_i$, then $\mathscr{F}$ coincides with the Lagrangian function $L(t, x, \dot{x})$ evaluated along the trajectories of the system.

\begin{theo}\label{sdfgtty}Let $\gamma \subset \tilde N$ be a closed curve. Then the integral
$$
\int_{g^\tau(\gamma)} \alpha
$$
is independent of $\tau$.
\end{theo}
This fact is a direct consequence of Theorem \ref{sdgffrt}.

Remark \ref{xfgb660} implies that the form $d\alpha$ is an integral invariant: $L_{\tilde w}d\alpha=0.$
\begin{theo}\label{dxfgyi}Let $\Sigma \subset \tilde N$ be an admissible two-dimensional surface. Then the integral
$$
\int_{g^\tau(\Sigma)} d\alpha
$$
is independent of $\tau$.\end{theo}
This fact is a direct consequence from Theorem \ref{sdfr5678}.
\begin{theo}\label{sdagft678}
Let $\Sigma \subset N$ be a two-dimensional surface. Then the integral
$$
\int_{G^t_{t_0}(\Sigma)} \beta
$$
is independent of $t$ and $t_0$.\end{theo}
\subsubsection*{Proof} This fact is a direct consequence of Theorem \ref{zxcvdv001}.
Indeed, since $d_z\beta=0$ and by virtue of equation (\ref{sdgf5txsw3}), we obtain:
$$
\frac{\partial \beta}{\partial t} + L_w \beta = d_z i_w \beta + i_w d_z \beta =  0.
$$
The theorem is proved.

\begin{cor}\label{zxcv567y}From theorem \ref {sdagft678} it follows that  $G^{t}_{t_0}$ preserves the form $\beta:$
\begin{equation}\label{dgf56y}(G^{t}_{t_0})_*\beta=\beta.\end{equation}
In other words, the shift along the trajectories of a Hamiltonian system is a \textbf{symplectic map}.

Moreover, equation (\ref{dgf56y}) implies that this shift preserves the volume in the phase space:
\[
(G^{t}_{t_0})_* \varpi = \varpi, \quad \varpi = \underbrace{\beta \wedge \ldots \wedge \beta}_{m \text{ times}}.
\]
\end{cor}

\begin{theo}If $\gamma \subset N$ is a closed curve, then the integral
$$
\int_{G^t_{t_0}(\gamma)} \zeta, \quad \text{where} \quad \zeta = p_i dx^i,
$$
is independent of $t$ and $t_0$.\end{theo}
This fact is a direct consequence of Theorem \ref{zxcvdv001}. Indeed, using the Cartan formula, we have:
$$
\frac{\partial \zeta}{\partial t} + L_w \zeta = d_z i_w \zeta + i_w d_z \zeta = d_z \left( p_i \frac{\partial H}{\partial p_i} - H \right).
$$

\begin{probl}\label{xcfvbgff}
Consider a Hamiltonian system with a Hamiltonian function
$H(t,x,p)$, where $x=(x^1,\ldots,x^m)$, $p=(p_1,\ldots,p_m)\in\mathbb{R}^m$, and $t\in\mathbb{R}$.
Assume this Hamiltonian to be positively homogeneous of degree one in the momenta:
$$
H(t,x,\lambda p) = \lambda H(t,x,p) \quad \forall \lambda > 0.
$$
Recall that such functions satisfy the Euler homogeneous function theorem:
$$
H = p_i \frac{\partial H}{\partial p_i}.
$$

Let $\gamma \subset N$ be an admissible one-dimensional manifold (not necessarily closed).
Show that the integral
$$
\int_{G_{t_0}^t(\gamma)} p_i dx^i
$$
is independent of $t$ and $t_0$.

\end{probl}

\begin{probl}\label{zsdvf55t}Let the Hamiltonian $H$ be as defined in Problem \ref{xcfvbgff}.
Fix the initial position $x(0) = \hat{x}$ and treat the initial momentum $p(0) = \hat{p}$ as a parameter.

Let the solution to the Hamiltonian equations be denoted by $x(t, \hat{p})$ and $p(t, \hat{p})$. Prove the equality:
\begin{equation}\label{vfrty}
\frac{\partial x^i(t, \hat{p})}{\partial \hat{p}_s} p_i(t, \hat{p}) = 0, \quad s = 1, \ldots, m.
\end{equation}

\textit{Hint:} Verify the equality for $t=0$ and, by direct calculation, show that
$$
\frac{d}{dt} \left( \frac{\partial x^i(t, \hat{p})}{\partial \hat{p}_s} p_i(t, \hat{p}) \right) = 0.
$$
\end{probl}

\section{Reduction of the Order of a Hamiltonian System by Use of the Energy Integral}Let us point out a consequence of Theorem \ref{cccfg}.

Assume that the Hamiltonian $H$ does not depend on $t$. Then $H$ is a first integral. Assume that the energy level
$$E_h = \{H(z) = h\} \subset \tilde{N}$$
is non-degenerate: $dH|_{E_h} \neq 0$.

Without loss of generality, we can assume that in any sufficiently small domain of the manifold $E_h$, the following inequality holds:
$$\frac{\partial H}{\partial p_1} \neq 0.$$
We can always achieve this by canonical permutations.\footnote{Changes of variables $(x, p) \mapsto (X, P)$ of the type
$$X^i = -p_i, \quad P_i = x^i, \quad P_s = p_s, \quad X^s = x^s, \quad s \neq i$$
or
$$X^i = x^j, \quad X^j = x^i, \quad P_i = p_j, \quad P_j = p_i$$
and compositions of such changes are called canonical permutations. From the results of Section \ref{sdgfkkkty}, it follows that canonical permutations preserve the Hamiltonian form of the equations.}

Then the Implicit Function Theorem implies that the surface $E_h$ is described locally by the graph
\begin{equation}\label{zsdvffghyu}
p_1 = g(x^1, \dots, x^m, p_2, \dots, p_m, h),
\end{equation}
where $t, x^1, \dots, x^m, p_2, \dots, p_m$ are local coordinates on $E_h$.

Therefore, we obtain:
$$\alpha|_{E_h} = \sum_{k=2}^m p_k dx^k - \big(-g(x^1, \dots, x^m, p_2, \dots, p_m, h)\big) dx^1 - h dt, \quad d(h dt) = 0.$$
Thus, the form $d(\alpha|_{E_h})$ is annihilated by the vector field associated with the following system:
\begin{align}
\frac{d t}{d T} &= b(t, x^1, \dots, x^m, p_1, \dots, p_m), \nonumber \\
\frac{d x^1}{d T} &= 1, \nonumber \\
\frac{d p_i}{d T} &= \frac{\partial g}{\partial x^i}, \quad \frac{d x^i}{d T} = -\frac{\partial g}{\partial p_i}, \quad i = 2, \dots, m \label{zxdff55bhy}
\end{align}
for an arbitrary function $b$.

Equations \eqref{zxdff55bhy} have Hamiltonian form. Their trajectories, together with the relation $x^1 = T + T_0$ and equation \eqref{zsdvffghyu}, are the projections of the trajectories of system \eqref{sxdfgyui} from $E_h$ onto the phase space $N$.

\section{Characteristic Property of the Hamilton-Jacobi Equation}
\begin{theo}\label{zdvfccc}

\begin{enumerate}
    \item Assume that a function $S = S(t, x)$ satisfies the Hamilton-Jacobi equation:
    \[
    H\left(t, x, \frac{\partial S}{\partial x}\right) + \frac{\partial S}{\partial t} = 0.
    \]
    Then the graph
    \[
    \Gamma = \left\{ p_i = \frac{\partial S}{\partial x^i}(t, x), \quad i = 1, \dots, m \right\} \subset \tilde{N}
    \]
    is an $(m+1)$-dimensional invariant manifold of system (\ref{sxdfgyui}).

    \item Conversely, let a graph
    \[
    \Gamma = \left\{ p = \frac{\partial S}{\partial x} \right\} \subset \tilde{N}, \quad S = S(t, x)
    \]
    be an invariant surface of system (\ref{sxdfgyui}). Then there exists a function $\psi = \psi(t)$ such that $S$ satisfies the equation
    \[
    H\left(t, x, \frac{\partial S}{\partial x}\right) + \frac{\partial S}{\partial t} = \psi(t).
    \]
    (Consequently, the function $\tilde{S} = S - \int \psi(t) \, dt$ satisfies the Hamilton-Jacobi equation.)

    \item If $S$ is a solution to the Hamilton-Jacobi equation, then $\alpha|_{\Gamma} = dS$.
\end{enumerate}
\end{theo}
\subsubsection*{Proof}The last assertion is trivial. Let us prove item 1).

Observe that $(t, x)$ are local coordinates on the manifold $\Gamma$. Let $x(t)$ be a solution to the system
\[
\dot{x}^i = \frac{\partial H}{\partial p_i} \left( t, x, \frac{\partial S}{\partial x}(t, x) \right).
\]
We shall show that
\[
x(t), \quad p_i(t) = \frac{\partial S}{\partial x^i}(t, x(t))
\]
is a solution to the Hamilton equations.

Indeed, let us introduce the function
\[
F(t, x) = H \left( t, x, \frac{\partial S}{\partial x}(t, x) \right).
\]
Then we have
\begin{equation} \label{xzcv3w4et}
\frac{\partial F}{\partial x^i} = \frac{\partial H}{\partial x^i} + \frac{\partial H}{\partial p_s} \frac{\partial^2 S}{\partial x^i \partial x^s}.
\end{equation}
Now we calculate $\dot{p}_k$:
\begin{align}
\dot{p}_k &= \frac{d}{dt} \left( \frac{\partial S}{\partial x^k}(t, x(t)) \right) = \frac{\partial^2 S}{\partial t \partial x^k} + \frac{\partial^2 S}{\partial x^k \partial x^s} \dot{x}^s \nonumber \\
&= \frac{\partial^2 S}{\partial t \partial x^k} + \frac{\partial^2 S}{\partial x^k \partial x^s} \frac{\partial H}{\partial p_s} = \frac{\partial^2 S}{\partial t \partial x^k} + \left( \frac{\partial F}{\partial x^k} - \frac{\partial H}{\partial x^k} \right) \nonumber \\
&= \frac{\partial}{\partial x^k} \left( \frac{\partial S}{\partial t} + F \right) - \frac{\partial H}{\partial x^k} = -\frac{\partial H}{\partial x^k}. \nonumber
\end{align}

To prove item 2), take a solution $(x(t), p(t))$ to the Hamilton equations and differentiate the equality
\[
p_i(t) = \frac{\partial S}{\partial x^i}(t, x(t))
\]
with respect to $t$:
\[
-\frac{\partial H}{\partial x^i} = \frac{\partial^2 S}{\partial x^i \partial x^r} \frac{\partial H}{\partial p_r} + \frac{\partial^2 S}{\partial x^i \partial t}.
\]
Using formula (\ref{xzcv3w4et}), we obtain
\[
\frac{\partial}{\partial x^i} \left( \frac{\partial S}{\partial t} + F \right) = 0.
\]
The theorem is proved.

Item 3) implies two consequences.
Firstly, we have $d\alpha|_{\Gamma} = 0$.
Secondly, if $(\tau, x(\tau), p(\tau))$ is a trajectory of (\ref{sxdfgyui}) that belongs to the manifold $\Gamma$, then the following equality holds:
\begin{align}
S(t, x(t)) - S(t_0, x(t_0)) &= \int_{t_0}^t \left( p_i(\tau) \dot{x}^i(\tau) - H(\tau, x(\tau), p(\tau)) \right) d\tau \nonumber \\
&= \int_{t_0}^t L(\tau, x(\tau), \dot{x}(\tau)) \, d\tau, \label{dfg66}
\end{align}
where $L$ is the Lagrangian.

The second equality holds only if the Legendre transform of $H$ is well-defined. Indeed,
\[
p_i(\tau) \dot{x}^i(\tau) - H(\tau, x(\tau), p(\tau)) = p_i(\tau) \frac{\partial H}{\partial p_i} - H(\tau, x(\tau), p(\tau)).
\]
For systems in classical mechanics, there are no issues with the Legendre transform.

Note also that if the Hamiltonian $H$ is a homogeneous function in $p$ (as in Problem \ref{xcfvbgff}), then the first line of equality (\ref{dfg66}) takes the form:
\[
S(t, x(t)) = S(t_0, x(t_0)).
\]

\section{Eikonal Equation and Gauss's Lemma}
\subsection{Eikonal Equation}
Let $M$ be a Riemannian manifold with local coordinates $x = (x^1, \dots, x^m)$, and let $g_{ij}(x)$ denote the components of the metric tensor.

Let $x(t)$ be a solution to the Euler-Lagrange equations with the Lagrangian
\begin{equation}
\label{dfgnn}
L = \frac{1}{2} g_{ij}(x) \dot{x}^i \dot{x}^j, \quad H = \frac{1}{2} g^{ij}(x) p_i p_j,
\end{equation}
where the initial conditions and conjugate momenta are given by
\[
x(0) = \hat{x}, \quad p_i = \frac{\partial L}{\partial \dot{x}^i} = g_{ij}(x) \dot{x}^j.
\]

Assume that the solution $x(t)$ has energy $1/2$:
\begin{equation}
\label{cfr4567}
|\dot{x}(t)|^2 = g_{ij}(x(t)) \dot{x}^i(t) \dot{x}^j(t) = 1.
\end{equation}
Under this assumption, the curve $x = x(t)$ defines a geodesic, and $t$ is the arc-length parameter.

\begin{theo}\label{sxdgffty}Suppose a function $f \colon M \to \mathbb{R}$ satisfies the eikonal equation
\begin{equation}
\label{dsfghjkl}
|\nabla f|^2 = g^{ij} \frac{\partial f}{\partial x^i} \frac{\partial f}{\partial x^j} = 1,
\end{equation}
where the gradient is given by $\nabla f = g^{ij} \frac{\partial f}{\partial x^i} \frac{\partial}{\partial x^j}$. Let $x(t)$ be a geodesic issued from a level surface
\[
\Psi_{\hat{x}} = \{x \in M \mid f(x) = f(\hat{x})\}, \quad \hat{x} \in \Psi_{\hat{x}},
\]
such that it is orthogonal to this surface at the initial point:
\begin{equation}
\label{sdfgg66}
\frac{\partial f}{\partial x^i}(\hat{x}) = g_{ij}(\hat{x}) \dot{x}^j(0).
\end{equation}
Then the following properties hold:
\begin{enumerate}
    \item The geodesic remains orthogonal to each level surface $\Psi_{x(t)}$ it intersects:
    \[
    \frac{\partial f}{\partial x^i}(x(t)) = g_{ij}(x(t)) \dot{x}^j(t);
    \]
    \item The value of the function along the geodesic satisfies:
    \[
    f(x(t)) - f(\hat{x}) = t.
    \]
\end{enumerate}
\end{theo}
\subsubsection*{Proof of Theorem \ref{sxdgffty}}The function $S(t,x) = f(x) - t/2$ satisfies the Hamilton--Jacobi equation with the Hamiltonian given in \eqref{dfgnn}. Moreover, due to condition \eqref{sdfgg66} and Theorem~\ref{zdvfccc}, the geodesic $x(t)$ belongs to the Lagrangian manifold $\Gamma$ defined by:
\[
p_i(t) = g_{ij}(x(t)) \dot{x}^j(t) = \frac{\partial S}{\partial x^i}(t, x(t)) = \frac{\partial f}{\partial x^i}(x(t)).
\]
This implies that the velocity vector $\dot{x}(t)$ is orthogonal to the level surface $\Psi_{x(t)}$.

By virtue of formula \eqref{cfr4567}, we have
\[
\int_0^t L \, d\tau = \frac{t}{2}.
\]
Consequently, from the relation for the action (or formula \eqref{dfg66}), we find
\[
S(t, x(t)) - S(0, \hat{x}) = \frac{t}{2}.
\]
Substituting $S(t,x) = f(x) - t/2$ into this equality, we obtain
\[
\left(f(x(t)) - \frac{t}{2}\right) - f(\hat{x}) = \frac{t}{2},
\]
which directly implies item~2) of the theorem.

Theorem~\ref{sxdgffty} is proved.

\begin{probl}\label{sdggg0}

Suppose a function $f \colon M \to \mathbb{R}$ satisfies the eikonal equation \eqref{dsfghjkl}. Show that all solutions to the autonomous system of ordinary differential equations
\begin{equation}
\label{xvbbb}
\dot{x} = \nabla f(x)
\end{equation}
are geodesics.
\end{probl}

\subsubsection*{Solution to Problem \ref{sdggg0}}The proof follows the previous argument almost literally.

We utilize the same solution to the Hamilton--Jacobi equation: $S(t,x) = f(x) - t/2$. For equation \eqref{xvbbb}, we fix an arbitrary initial value $x(0) = \hat{x}$ and set the initial momentum to $\hat{p} = \frac{\partial S}{\partial x}(0,\hat{x})$.

Let $(x, p)(t)$ be a solution to the Hamilton equations with the initial conditions $(\hat{x}, \hat{p})$. This implies, in particular, that $x(t)$ is a solution to the Euler--Lagrange equations; hence, $x(t)$ is a geodesic. By Theorem~\ref{zdvfccc}, we have
\[
p_i(t) = g_{ij}(x(t)) \dot{x}^j(t) = \frac{\partial S}{\partial x^i}(t, x(t)) = \frac{\partial f}{\partial x^i}(x(t)).
\]
Consequently, $x(t)$ is a solution to \eqref{xvbbb}.

There is, however, a more concise way to solve Problem~\ref{sdggg0}. Any solution to \eqref{xvbbb} minimizes the following action functional:
\[
I[x] = \int_{t_1}^{t_2} \frac{1}{2} |\dot{x} - \nabla f|^2 \, dt.
\]
Note that the integrand can be expanded as
\[
\frac{1}{2} |\dot{x} - \nabla f|^2 = \frac{1}{2} |\dot{x}|^2 + \frac{1}{2} |\nabla f|^2 - \langle \dot{x}, \nabla f \rangle = \frac{1}{2} |\dot{x}|^2 + \frac{1}{2} - \frac{df}{dt}.
\]
Since $\frac{1}{2}$ and the total derivative $df/dt$ do not affect the Euler--Lagrange equations, any solution to \eqref{xvbbb} must also be a solution to the equations of motion for the Lagrangian $L = \frac{1}{2} |\dot{x}|^2$.

\begin{probl}Let $f(x, a)$ be an $\ell$-parametric family of solutions to the eikonal equation \eqref{dsfghjkl}, where $a = (a^1, \dots, a^\ell) \in \mathbb{R}^\ell$ denotes the parameters.

Assume that the system of equations
\begin{equation}
\label{envelope_cond}
\frac{\partial f}{\partial a}(x, a) = 0
\end{equation}
possesses a smooth solution $a = \theta(x)$.

Show that the envelope defined by
\[
\tilde{f}(x) := f(x, \theta(x))
\]
is also a solution to equation \eqref{dsfghjkl} (see \cite{ev}).

Furthermore, show that for any given value of $a$, the level surface $\{\tilde{f}(x) = c\}$ is tangent to the corresponding surface $\{f(x, a) = c\}$ of the family at their points of intersection.
\end{probl}

\subsection{Gauss's Lemma}Consider the result of Problem~\ref{zsdvf55t} from a Riemannian geometry perspective.

The Hamiltonian
$$
\mathcal{H} = \sqrt{g^{ij}(x) p_i p_j}
$$
is a homogeneous function of degree one in the momenta. The function $\mathcal{H}$ is a first integral of the system. Let $(x(t), p(t))$ be a solution with a fixed initial position $x(0) = \hat{x}$ and satisfying the condition
\begin{equation}
\label{xcvbvvgf}
\mathcal{H} = 1.
\end{equation}
We regard the initial momentum $p(0) = \hat{p}$ as a parameter.

It is straightforward to see that such solutions also satisfy the system with the Hamiltonian given in \eqref{dfgnn}. In this case, condition \eqref{xcvbvvgf} is equivalent to \eqref{cfr4567}. In particular, $t$ is the arc-length parameter, and solutions $x(t)$ of the type described above are geodesics on $M$ emanating from the point $\hat{x}$.

Consider a sphere of radius $r > 0$ in the tangent space $T_{\hat{x}}M$:
\[
C_r = \{v = (v^1, \dots, v^m) \in T_{\hat{x}}M \mid |v|^2 = g_{ij}(\hat{x}) v^i v^j = r^2\}.
\]
Each vector $v \in C_r$ defines a unique geodesic $x_v(t)$ satisfying the initial conditions $x_v(0) = \hat{x}$ and $\dot{x}_v(0) = v$.

Let $r$ be chosen sufficiently small so that any two distinct geodesics from the family $\{x_v\}_{v \in C_r}$ do not intersect for $t \in (0, 1]$. Under this condition, the set
\[
\Sigma = \{x_v(1) \in M \mid v \in C_r\}
\]
is a smooth submanifold of $M$ of dimension $m-1$.
\begin{theo}The velocity vector $\dot{x}_v(1) \in T_{x_v(1)}M$, where $|v|=r$, is orthogonal to the tangent space of the submanifold $\Sigma$ at the point $x_v(1)$.\end{theo}
Indeed, consider an arbitrary smooth curve $x(\xi) = x_{v(\xi)}(1)$ on $\Sigma$, where $|v(\xi)| = r$ and $\xi \in \mathbb{R}$ is a parameter along the curve.

A tangent vector to $\Sigma$ at the point $x_v(1)$ is given by:
\[
\frac{d x^i_{v(\xi)}(1)}{d\xi} = \frac{\partial x^i_v}{\partial v^r} \frac{dv^r}{d\xi} = \frac{\partial x^i_v}{\partial \hat{p}_j} \frac{\partial \hat{p}_j}{\partial v^r} \frac{dv^r}{d\xi} = \frac{\partial x^i_v}{\partial \hat{p}_j} g_{jr}(\hat{x}) \frac{dv^r}{d\xi} = \frac{\partial x^i_v}{\partial \hat{p}_j} \frac{d\hat{p}_j}{d\xi}.
\]

By virtue of \eqref{vfrty}, the following orthogonality condition holds:
\[
g_{il}(x_v(1)) \dot{x}^l_v(1) \frac{\partial x^i_v}{\partial \hat{p}_j} = 0.
\]
Consequently, it follows that
\[
g_{il}(x_v(1)) \dot{x}^l_v(1) \frac{d}{d\xi} x^i_v(1) = 0,
\]
which means the velocity vector $\dot{x}_v(1)$ is orthogonal to every tangent vector of $\Sigma$. This completes the proof of the theorem.

\section{Canonical Transformations. Generating Functions}\label{sdgfkkkty}
\subsection{Canonical Transformations }Let $f = f(t, x, p)$ be a function
defined on the extended phase space $\tilde{N}$ of the system \eqref{sdfggg}. Recall the notation:
\[
d_z f = \frac{\partial f}{\partial x^i} dx^i + \frac{\partial f}{\partial p_i} dp_i, \quad df = d_z f + \frac{\partial f}{\partial t} dt.
\]
\begin{df}[\cite{arn}]A transformation\footnote{We use the word ``transformation'' as a synonym for the phrase ``change of variables.''} of the manifold $\tilde{N}$
\[
(t, x, p) \mapsto (t, X, P), \quad X = X(t, x, p), \quad P = P(t, x, p)
\]
is said to be \textit{canonical} if
\begin{equation}
\label{xcfbghoo}
d_z P_i \wedge d_z X^i = dp_i \wedge dx^i.
\end{equation}
\end{df}
In the following, we treat $t, x, p$ as independent variables, such that
\[
d_z p_i = dp_i, \quad d_z x^i = dx^i, \quad d_z t = 0.
\]

Loosely speaking, a canonical transformation is a family of symplectic maps parameterized by $t$. For instance, according to Corollary~\ref{zxcv567y}, the mapping $(t, x, p) \mapsto (t, G^t_{t_0}(x, p))$ is a canonical transformation.

Equation \eqref{xcfbghoo} is locally equivalent to the existence of a function $S(t, x, p)$ such that
\begin{equation}
\label{xcvbgf56yhn}
p_i dx^i - P_i d_z X^i = d_z S(t, x, p).
\end{equation}
Indeed, this follows from the fact that the exterior derivative vanishes: $d_z(p_i dx^i - P_i d_z X^i) = 0$.

Equality \eqref{xcvbgf56yhn} can be rewritten in the form
\begin{equation}
\label{sdfsdfv56z}
-P_i dX^i + P_i \frac{\partial X^i}{\partial t} dt + p_i dx^i = dS - \frac{\partial S}{\partial t} dt.
\end{equation}

\begin{theo}\label{zdsfsadfgll}In the new coordinates $(X, P)$, the Hamiltonian system \eqref{sdfggg} preserves its canonical form:
$$
\dot{P}_i = -\frac{\partial K}{\partial X^i}, \quad \dot{X}^i = \frac{\partial K}{\partial P_i},
$$
where the new Hamiltonian $K(t, X, P)$ is given by
$$
K(t, X, P) = \left( P_i \frac{\partial X^i}{\partial t} + \frac{\partial S}{\partial t} + H \right) \bigg|_{(x,p) \to (X,P)}.
$$
\end{theo}
Indeed, formula \eqref{sdfsdfv56z} implies that
\begin{equation}
\label{zxdvgtyyyy}
\alpha = p_i dx^i - H dt = P_i dX^i - K dt + dS,
\end{equation}
and the assertion of Theorem~\ref{zdsfsadfgll} follows directly from Theorem~\ref{cccfg}.
\begin{cor}If a canonical transformation does not depend on time,
\[
P = P(x, p), \quad X = X(x, p),
\]
then the Hamiltonian $H$ transforms as a scalar function on the manifold $N$:
\begin{equation}
\label{time_indep_K}
K(t, X, P) = H\big(t, x(X, P), p(X, P)\big).
\end{equation}\end{cor}

\subsection{Generating Functions}

\begin{df} A canonical transformation $(t, x, p) \mapsto (t, X, P)$ is said to be \textit{free} if
\begin{equation}
\label{sfffrt5}
\det \left( \frac{\partial X^i}{\partial p_j} \right) \neq 0.
\end{equation}\end{df}
In this case, by the Implicit Function Theorem, the variables $(t, x, X)$ serve as local coordinates on $\tilde{N}$. In particular, we can express the function $S$ as $S = S_1(t, x, X)$.

Formula \eqref{zxdvgtyyyy} then takes the form
\begin{equation}
\label{dth6878}
p_i dx^i - H dt = P_i dX^i - K dt + \frac{\partial S_1}{\partial t} dt + \frac{\partial S_1}{\partial x^i} dx^i + \frac{\partial S_1}{\partial X^i} dX^i,
\end{equation}
which implies the following relations:
\begin{align}
p_i &= \frac{\partial S_1}{\partial x^i}, \quad P_i = -\frac{\partial S_1}{\partial X^i}; \label{srfrrf} \\
K &= H + \frac{\partial S_1}{\partial t}. \label{dfff}
\end{align}
The condition \eqref{sfffrt5} can be written as
\begin{equation}
\label{4tt455}
\det \left( \frac{\partial^2 S_1}{\partial x^i \partial X^j} \right) \neq 0.
\end{equation}
Conversely, if the inequality \eqref{4tt455} holds, then the relations \eqref{srfrrf} define a canonical transformation.

Indeed, let us verify the non-degeneracy of the mapping $(x, p) \mapsto (X, P)$. We can represent this transformation as a composition of two steps:
\[
(x, p) \mapsto (x', X') \mapsto (X, P).
\]
The first step is the inverse of the transformation defined by:
\[
x = x', \quad p = \frac{\partial S_1}{\partial x}(t, x', X').
\]
The Jacobian of this change of variables is
\[
\det \left( \frac{\partial^2 S_1}{\partial x^i \partial X^j} \right)^{-1}.
\]
The second step is given by the formulas:
\[
X = X', \quad P = -\frac{\partial S_1}{\partial X}(t, x', X').
\]
Its Jacobian is
\[
\det \left( \frac{\partial^2 S_1}{\partial x^i \partial X^j} \right).
\]

Since the total Jacobian (the product of these two) is non-zero, the transformation is well-defined. Therefore, if the function $S_1$ satisfies \eqref{4tt455}, it defines a free transformation which is uniquely determined by the relations \eqref{srfrrf}.

The function $S_1$ is called a \textit{generating function}.

In the new variables $(X, P)$, the Hamiltonian $K$ is given by \eqref{dfff}.

\begin{rem}\label{xfgtrpd}Assume that we have a generating function $S_1(t, x, X)$ that satisfies the Hamilton--Jacobi equation:
$$
H\left( t, x, \frac{\partial S_1}{\partial x} \right) + \frac{\partial S_1}{\partial t} = 0.
$$
Then, the transformation $(t, x, p) \mapsto (t, X, P)$ is determined, and in the new variables $(X, P)$, the new Hamiltonian $K$ vanishes identically ($K \equiv 0$).

Consequently, the equations of motion in the new coordinates take the form $\dot{X} = 0$ and $\dot{P} = 0$, meaning the system is integrated immediately.\end{rem}

Consider a composition of canonical transformations:
\[
(x, p) \mapsto (\tilde{X}, \tilde{P}) \mapsto (X, P).
\]
The first transformation is defined by a generating function $S_1(t, x, \tilde{X})$ via the relations:
\begin{equation}
\label{dfh668}
p_i = \frac{\partial S_1}{\partial x^i}, \quad \tilde{P}_i = -\frac{\partial S_1}{\partial \tilde{X}^i}.
\end{equation}
The second transformation is a canonical permutation:
\[
\tilde{X}^i = P_i, \quad \tilde{P}_i = -X^i.
\]
Thus, the change $(x, p) \mapsto (X, P)$ is expressed by means of a generating function $S_2(t, x, P) := S_1(t, x, P)$. Formulas \eqref{dfh668} then imply:
$$
p_i = \frac{\partial S_2}{\partial x^i}, \quad X^i = \frac{\partial S_2}{\partial P_i}.
$$
In this case, the non-degeneracy condition \eqref{4tt455} takes the form
$$
\det \left( \frac{\partial^2 S_2}{\partial x^i \partial P_j} \right) \neq 0,
$$
and the condition \eqref{sfffrt5} for a free transformation is written as
$$
\det \left( \frac{\partial P_i}{\partial p_j} \right) \neq 0.
$$

By applying canonical permutations to various subsets of the conjugate pairs $(\tilde{X}^i, \tilde{P}_i)$, we can obtain $2^m$ different types of generating functions.

\begin{rem} The identity mapping
\[
P = p, \quad X = x
\]
is obtained via the $S_2$-type generating function $S_2 = P_i x^i$.

\end{rem}

\begin{df}A function $S = S(t, x, b)$ depending on $m$ parameters $b = (b_1, \dots, b_m)$ is called a \textit{complete integral} of the Hamilton--Jacobi equation
$$
H\left(t, x, \frac{\partial S}{\partial x}\right) + \frac{\partial S}{\partial t} = 0,
$$
provided that it satisfies the equation for all admissible values of $b$ and the non-degeneracy condition
$$
\det \left( \frac{\partial^2 S}{\partial x^i \partial b_j} \right) \neq 0
$$
holds.\end{df}Therefore, if a complete integral is provided, the system of Hamilton's equations can be integrated in closed form. See Remark~\ref{xfgtrpd}.

The vector $b$ can be interpreted as a set of new coordinates, in which case we obtain the generating function $S_1$. Alternatively, if we treat $b$ as a set of new momenta, we obtain $S_2$. More generally, the vector $b$ can consist of a mixture of both new coordinates and new momenta.

Currently, most known integrable Hamiltonian systems are integrated using the \textit{separation of variables} procedure. For the autonomous case ($H = H(z)$), this means that in certain suitable coordinates, the Hamilton--Jacobi equation
$$
H\left( x, \frac{\partial S}{\partial x} \right) = K(b)
$$
admits a complete integral of the form
$$
S(x, b) = \sum_{k=1}^m S_k(x^k, b).
$$

\section{Hamiltonian Vector Field Straightening Theorem}
Assume that the system \eqref{sdfggg} is autonomous, i.e., $H = H(z)$.
\begin{theo}\label{zxdbfr5}
Assume that the Hamiltonian $H$ is non-degenerate at a point $\tilde{z} \in N$:
$$
dH(\tilde{z}) \neq 0.
$$
Then, in some neighborhood $U$ of the point $\tilde{z}$, there exist canonical coordinates
\[
Z = (X, P) = (X^1, \dots, X^m, P_1, \dots, P_m), \quad Z = Z(x, p)
\]
such that:
\begin{enumerate}
    \item $dp_i \wedge dx^i = dP_i \wedge dX^i$ (the transformation is canonical);
    \item In the coordinates $Z$, the Hamiltonian takes the form $H = X^1$.
\end{enumerate}\end{theo}
\subsubsection*{Proof} Without loss of generality, let $\tilde{z}=0$ and
\begin{equation} \label{xfg56f}
H(0)=0, \quad \frac{\partial H}{\partial p_1}(0) \neq 0.
\end{equation}
Then, by the Implicit Function Theorem, the equation
\begin{equation}\label{xfg6000} H(x,p) = X^1 \end{equation}
has a solution
\[ p_1 = \phi(x, p_2, \dots, p_m, X^1), \quad H(x, \phi, p_2, \dots, p_m) = X^1 \]
for sufficiently small $|X^1|$, $|x|$, and $|p|$ (cf. \eqref{zsdvffghyu}).
From \eqref{xfg56f}, \ref{xfg6000} we obtain
\begin{equation} \label{zsdg5tl}
\phi(0)=0, \quad \frac{\partial H}{\partial p_1}(0) \cdot \frac{\partial \phi}{\partial X^1}(0) = 1.
\end{equation}

Consider the following Cauchy problem:
\begin{equation} \label{zdfff}
\frac{\partial S}{\partial x^1} = \phi\left(x, \frac{\partial S}{\partial x^2}, \dots, \frac{\partial S}{\partial x^m}, X^1\right), \quad S\big|_{x^1=0} = \sum_{k=2}^m x^k X^k, \quad S = S(x, X).
\end{equation}
As shown in Section~\ref{zdvgfsdegf5poji}, this Cauchy problem has a solution $S = S(x, X)$ for sufficiently small $|x|$ and $|X|$.

Thus, equation \eqref{zdfff} is equivalent to the equation
\[ H\left(x, \frac{\partial S}{\partial x}\right) = X^1. \]
To complete the proof, it remains to verify that $S(x, X)$ is a generating function of a transformation $(x, p) \mapsto (X, P)$:
\[ p = \frac{\partial S}{\partial x}, \quad P = -\frac{\partial S}{\partial X}. \]
The values $x=0, p=0$ correspond to $X=0, P=0$.

By virtue of \eqref{zdfff}, we have
\[
\frac{\partial^2 S}{\partial x \partial X} \bigg|_{x=X=0} =
\begin{bmatrix}
\frac{\partial \phi}{\partial X^1} & \frac{\partial^2 S}{\partial x^1 \partial X^2} & \dots & \frac{\partial^2 S}{\partial x^1 \partial X^m} \\
0 & 1 & \dots & 0 \\
\vdots & \vdots & \ddots & \vdots \\
0 & 0 & \dots & 1
\end{bmatrix}, \quad
\det \frac{\partial^2 S}{\partial x \partial X} \bigg|_{x=X=0} = \frac{\partial \phi}{\partial X^1}.
\]
From \eqref{zsdg5tl}, it follows that
\[ \frac{\partial \phi}{\partial X^1}(0) \neq 0. \]

The theorem is proved.

\section{The Poincare Section on the Energy Level, \cite{TrZ} }
Recall that the Hamiltonian vector field $w$ is defined by formula \eqref{xdfsg500i}. We assume that $H$ does not depend on $t$.

Let $E_h = \{z \in N \mid H(z) = h\}$ denote an energy level, assuming $dH|_{E_h} \neq 0$.

Let $Y \subset E_h$ be a hypersurface of dimension $\dim Y = 2m - 2$ such that $w(z) \notin T_z Y$ for all $z \in Y$. Consider the restriction $\beta' = \beta|_Y$.

\begin{theo}\label{xcvf5tyutyt7}
The form $\beta'$ is nondegenerate; thus, $(Y, \beta')$ is a symplectic manifold.
\end{theo}

\subsubsection*{Proof}
Take a point $z \in Y$. In a neighborhood of this point, introduce canonical coordinates in accordance with Theorem \ref{zxdbfr5}. In these coordinates, the integral curves are given by:
$$
P_1(t) = -t + P_1(0), \quad P_j(t) = P_j(0), \quad X^i(t) = X^i(0), \quad j=2, \dots, m,
$$
and the energy level is $E_h = \{X^1 = h\}$.

Thus, $(X^2, \dots, X^m, P_1, \dots, P_m)$ serve as local coordinates on $E_h$, and the manifold $Y$ can be represented as a graph $P_1 = y(X^2, \dots, X^m, P_2, \dots, P_m)$. This implies that $(X^2, \dots, X^m, P_2, \dots, P_m)$ are local coordinates on $Y$.

The proof is completed by the following computation:
$$
\beta|_Y = \left. \sum_{i=1}^m dP_i \wedge dX^i \right|_Y = \sum_{k=2}^m dP_k \wedge dX^k.
$$
The theorem is proved.

Let $G^t: N \to N$ denote the phase flow of \eqref{sdfggg}. Recall that $G^t_* \beta = \beta$ (see formula (\ref{dgf56y})). Assume that the manifolds $Y_1, Y_2 \subset E_h$ satisfy the same properties as $Y$. As proved above, the pairs $(Y_i, \beta_i)$, where $\beta_i = \beta|_{Y_i}$, are symplectic manifolds.

Suppose that a trajectory starting at $z \in Y_1$ reaches the surface $Y_2$ at time $t = \tau(z) > 0$, such that:
$$
Q(z) := G^{\tau(z)}(z) \in Y_2.$$
\begin{theo}\label{sxdfg5jhjhjhyu6}
The mapping $Q \colon Y_1 \to Y_2$ is symplectic; that is, $Q_* \beta_2 = \beta_1$.\end{theo}
\subsubsection*{Proof}
Take a two-dimensional compact manifold $\Sigma \subset Y_1$ with a smooth boundary curve $\partial \Sigma$ without self-intersections and consider a cylinder
$$C = \bigcup_{z \in \Sigma} \{G^t(z) \mid t \in [0, \tau(z)]\}.$$
Since the form $\beta$ is closed, we have
$$\int_{\partial C} \beta = 0.$$
On the other hand, it follows that
$$\partial C = \Sigma \cup Q(\Sigma) \cup \sigma, \quad \sigma = \bigcup_{z \in \partial \Sigma} \{G^t(z) \mid t \in (0, \tau(z))\}.$$
By a suitable choice of orientations, we obtain
$$\int_\Sigma \beta = \int_{Q(\Sigma)} \beta + \int_\sigma \beta.$$
Since $\Sigma$ is an arbitrary compact manifold, to complete the proof, it is enough to show that
$$\beta|_\sigma = 0.$$
A basis in $T_z\sigma$ can be chosen as follows: $w, v$, where $w$ is the Hamiltonian vector field and $v$ is some other vector.

Let us calculate $\beta$ on the basis vectors by using formula (\ref{sdgf5txsw3}):
$$\beta(w, v) = -dH(v) = 0.$$
The last equality holds because $v \in T_z\sigma \subset T_z E_h$.

The theorem is proved.

The result of Theorem \ref{sxdfg5jhjhjhyu6} remains valid if $Y_1 = Y_2$.

\section{The Hamilton--Jacobi Equation in the General Case: The Method of Characteristics}
 \label{zdvgfsdegf5poji}
In the theory of PDEs, the equation
\begin{equation} \label{dthhh87}
u_t + f\left(t, x, u, \frac{\partial u}{\partial x}\right) = 0
\end{equation}
is also called the Hamilton--Jacobi equation.

Here, $f = f(t, x, \xi, p)$ is a scalar function of the variables
\[ t, \xi \in \mathbb{R}, \quad p = (p_1, \dots, p_m), \quad x = (x^1, \dots, x^m) \in \mathbb{R}^m. \]
Consider the following system of ODEs:
\begin{equation} \label{dh678}
\begin{aligned}
\dot{\xi} &= p_i \frac{\partial f}{\partial p_i} - f, \\
\dot{p}_i &= -\frac{\partial f}{\partial x^i} - \frac{\partial f}{\partial \xi} p_i, \\
\dot{x}^i &= \frac{\partial f}{\partial p_i}.
\end{aligned}
\end{equation}
The extended phase space $V$ of this system is the space of variables
\[ V = \{ (t, x, \xi, p) \in \mathbb{R}^{2m+2} \}. \]
If the function $f$ does not depend on $\xi$, then the equations for $x$ and $p$ in \eqref{dh678} decouple and form a system of Hamiltonian equations, while equation \eqref{dthhh87} becomes a Hamilton--Jacobi equation in the sense discussed above.

System \eqref{dh678} is called the characteristic system, and its solutions are called characteristics.

The proof of the following theorem is analogous to the proof of Theorem~\ref{zdvfccc}.
\begin{theo}\label{11122wwddf} Let $u=u(t,x)$ be a solution to \eqref{dthhh87}. Then the $(m+1)$-dimensional manifold (graph)
$$
G = \left\{ \xi=u(t,x), \quad p_i = \frac{\partial u}{\partial x^i}(t,x), \quad i=1,\dots,m \right\} \subset V
$$
is an invariant manifold of system \eqref{dh678}, and the following formula holds:
$$
(p_i dx^i - f dt - d\xi)\big|_{G} = 0,
$$
or equivalently,
$$
(p_i dx^i - f dt)\big|_{G} = du.
$$
\end{theo}This theorem provides a method for solving the Cauchy problem for equation \eqref{dthhh87}. Indeed, consider this equation with the following initial data:
\begin{equation} \label{dth56}
u\big|_{t=0} = \hat{u}(x).
\end{equation}

Define the initial conditions for system \eqref{dh678} as follows:
\[ x\big|_{t=0} = \hat{x}, \quad \xi\big|_{t=0} = \hat{u}(\hat{x}), \quad p\big|_{t=0} = \frac{\partial \hat{u}}{\partial x}(\hat{x}), \]
and let $(x, \xi, p) = (X, \Xi, P)(t, \hat{x})$ be the corresponding solution to system \eqref{dh678}.

From Theorem \ref{11122wwddf}, the solution $u(t, x)$ to the Cauchy problem \eqref{dthhh87}, \eqref{dth56} satisfies
$$
    \Xi(t,\hat x) = u(t, X(t,\hat x)).
$$
Thus, to determine $u(t,x)$, we must solve the equation $x = X(t,\hat x)$ for $\hat x$ and substitute the result into $\Xi(t,\hat x)$. This is possible for small $|t|$ due to the Implicit Function Theorem.

This approach to solving the Cauchy problem \eqref{dthhh87}, \eqref{dth56} is known as the method of characteristics.

\section{Appendix}
\subsection{Proof of Formula (\ref{zsdg00v})}
Let us recall this formula first:
$$di_v+i_vd=L_v.$$
First, we verify the formula for each point $x \in M$ such that $v(x) \neq 0$. It is a well-known result  \cite{taylor} that if $v(\tilde{x}) \neq 0$, then in some neighborhood of $\tilde{x}$ there exist local coordinates $x = (x^1, \dots, x^m)$ such that the vector field $v$ is represented as $v = (1,0,\ldots,0)$. The corresponding flow is given by:
\begin{equation}
\label{dgfhhh}
g^t(x) = (x^1 + t, x^2, \dots, x^m).
\end{equation}

By the linearity of $L_v, d,$ and $i_v$, it is sufficient to check the Cartan formula for monomials of the following two types:
\begin{enumerate}
    \item $\omega = a(x) dx^1 \wedge dx^{j_1} \wedge \dots \wedge dx^{j_{k-1}}, \quad 1 < j_1 < \dots < j_{k-1} \leq m$;
    \item $\gamma = b(x) dx^{l_1} \wedge \dots \wedge dx^{l_k}, \quad 1 < l_1 < \dots < l_k \leq m$.
\end{enumerate}

Consider Case 1. Case 2 is analogous but simpler. By direct calculation, we obtain:
\begin{equation*}
i_v \omega = a \, dx^{j_1} \wedge \dots \wedge dx^{j_{k-1}}, \quad d i_v \omega = \sum_{s=1}^m \frac{\partial a}{\partial x^s} dx^s \wedge dx^{j_1} \wedge \dots \wedge dx^{j_{k-1}}.
\end{equation*}
Furthermore,
\begin{align*}
d \omega &= \sum_{r=2}^m \frac{\partial a}{\partial x^r} dx^r \wedge dx^1 \wedge dx^{j_1} \wedge \dots \wedge dx^{j_{k-1}}, \\
i_v d \omega &= -\sum_{r=2}^m \frac{\partial a}{\partial x^r} dx^r \wedge dx^{j_1} \wedge \dots \wedge dx^{j_{k-1}}.
\end{align*}
Summing these results, the terms in the summation for $s \geq 2$ cancel out, leaving:
\begin{equation*}
d i_v \omega + i_v d \omega = \frac{\partial a}{\partial x^1} dx^1 \wedge dx^{j_1} \wedge \dots \wedge dx^{j_{k-1}}.
\end{equation*}

On the other hand, from the flow formula \eqref{dgfhhh}, it follows that:
\begin{align*}
L_v \omega &= \frac{d}{dt} \bigg|_{t=0} a(x^1+t, x^2, \dots, x^m) d(x^1+t) \wedge dx^{j_1} \wedge \dots \wedge dx^{j_{k-1}} \\
&= \frac{\partial a}{\partial x^1} dx^1 \wedge dx^{j_1} \wedge \dots \wedge dx^{j_{k-1}}.
\end{align*}
This proves the Cartan formula at each point of the set $F = \{x \in M \mid v(x) \neq 0\}$.

Since $F$ is an open set, and all involved functions are continuous, the Cartan formula remains valid on the closure $\overline{F}$. On the open set $N = M \setminus \overline{F}$, the vector field vanishes identically ($v|_N \equiv 0$). Consequently, both sides of the formula vanish, reducing to the trivial identity $0 = 0$. The proof is complete.

\subsection{Proof of Formula (\ref{Psfgb66})}
Let us restate the formula here for convenience:
$$
L_u i_v - i_v L_u = i_{[u, v]}.
$$

Since the operators $L_v$ and $i_u$ are linear, it is sufficient to obtain formula \eqref{Psfgb66} for monomials
$$\omega = a(x) dx^{i_1} \wedge \dots \wedge dx^{i_k}, \quad 1 \le i_1 < \dots < i_k\le m.$$

 In a neighborhood of a point $\tilde{x}$ where $u(\tilde{x}) \neq 0$, there exist local coordinates such that $u = (1, 0, \dots, 0)$. We assume that these coordinates, denoted by $x = (x^1, \dots, x^m)$, have already been introduced. Let us verify formula \eqref{Psfgb66} in such a neighborhood.

By  Definition \ref{sdg000d}, we have
\begin{equation*}
    L_u\omega = \frac{\partial a}{\partial x^1} dx^{i_1} \wedge \dots \wedge dx^{i_k},
\end{equation*}
and from formula \eqref{sswe4rfv}, it follows that
\begin{equation*}
    i_v L_u \omega = \frac{\partial a}{\partial x^1} \sum_{s=1}^k (-1)^{s-1} v^{i_s} dx^{i_1} \wedge \dots \wedge \widehat{dx^{i_s}} \wedge \dots \wedge dx^{i_k}.
\end{equation*}
Here, the hat symbol $\widehat{\dots}$ denotes the omission of the corresponding term.

For the same reason, it follows that
\begin{align}
L_u i_v \omega &= L_u \bigg( a \sum_{s=1}^k (-1)^{s-1} v^{i_s} dx^{i_1} \wedge \dots \wedge \widehat{dx^{i_s}} \wedge \dots \wedge dx^{i_k} \bigg) \nonumber \\
&= \sum_{s=1}^k \frac{\partial (a v^{i_s})}{\partial x^1} (-1)^{s-1} dx^{i_1} \wedge \dots \wedge \widehat{dx^{i_s}} \wedge \dots \wedge dx^{i_k}. \nonumber
\end{align}
Subtracting one formula from the other, we obtain
\begin{equation*}
L_u i_v \omega - i_v L_u \omega = a \sum_{s=1}^k \frac{\partial v^{i_s}}{\partial x^1} (-1)^{s-1} dx^{i_1} \wedge \dots \wedge \widehat{dx^{i_s}} \wedge \dots \wedge dx^{i_k}.
\end{equation*}
It remains to observe that in these specific coordinates, the Lie bracket is given by
\begin{equation*}
    [u,v]^i = \frac{\partial v^i}{\partial x^1},
\end{equation*}
which implies that
\begin{equation*}
    i_{[u,v]}\omega = L_u i_v \omega - i_v L_u \omega.
\end{equation*}

Thus, we have proved formula \eqref{Psfgb66} on the set
\begin{equation*}
    U = \{x \in M \mid u(x) \neq 0\}.
\end{equation*}
This set is clearly open. By the continuity of the functions involved, formula \eqref{Psfgb66} remains valid on its closure $\overline{U}$. Therefore, it remains only to verify the formula on the open set $M \setminus \overline{U}$. In this set, $u \equiv 0$, and the formula reduces to the trivial identity $0 = 0$.

\begin{rem}
Let $u(x), v(x)$ be vector fields defined on the $m$-dimensional manifold $M$, and let $\omega$ be a $1$-form on $M$. Formulas (\ref{zsdg00v}) and (\ref{Psfgb66}) imply
\begin{equation}\label{sdfgh--1}
d\omega(u,v)=L_u\big(\omega(v)\big)-L_v\big(\omega(u)\big)-i_{[u,v]}\omega.
\end{equation}
\end{rem}

\begin{rem}
Let $v_1(x), \ldots, v_m(x)$ be vector fields defined on $M$. Suppose these vector fields are linearly independent at each point $x \in M$ and
$$[v_i, v_j] = c_{ij}^k(x) v_k, \quad c_{ij}^k = -c_{ji}^k.$$
Such an expansion is valid due to the linear independence of the vector fields. Thus, this formula serves to define the functions $c_{ij}^k(x)$. Let us introduce $1$-forms $\omega^1, \ldots, \omega^m$ such that $\omega^i(v_j) = \delta_j^i$.

Formula (\ref{sdfgh--1}) implies that
$$i_{[v_k, v_s]} \omega^r = -i_{v_s} i_{v_k} d\omega^r = c_{ks}^r.$$
This identity leads to the Maurer--Cartan formula
\begin{equation}\label{d---sfgyy}
-d\omega^r = \frac{1}{2} c^r_{ij} \omega^i \wedge \omega^j.
\end{equation}
To verify this, it is sufficient to check formula (\ref{d---sfgyy}) on pairs of basis vectors $v_k, v_s$.
\end{rem}

\subsection{ Poincare Lemma}Assume that a smooth $k$-form ($k\ge 1$)
$$\omega = \sum_{i_1 < \dots < i_k} \omega_{i_1 \dots i_k} (x) dx^{i_1} \wedge \dots \wedge dx^{i_k}$$
is defined in a neighbourhood of the origin in $\mathbb{R}^m$.

\begin{theo}[Poincaré Lemma]
Suppose that $\omega$ is closed ($d\omega = 0$). Then $\omega$ is exact in some open ball $B$ centered at the origin. That is, there exists a $(k-1)$-form $\Omega$ defined on $B$ such that $d\Omega = \omega$.
\end{theo}
The form $\Omega$ is called a primitive form.
\begin{proof}
Consider a vector field $v(x) = x$. The corresponding flow is $g^t(x) = e^t x$, and the pullback of $\omega$ is given by
\begin{equation}\label{dfh700}
g^t_* \omega = e^{kt} \sum_{i_1 < \dots < i_k}
\omega_{i_1 \dots i_k}(e^t x) dx^{i_1} \wedge \dots \wedge dx^{i_k}.
\end{equation}

As a special case of  formula (\ref{zdggty}), we have
$$\frac{d}{dt} g^t_* \omega = g^t_* L_v \omega.$$
Applying Cartan's formula and the properties of the pullback, we obtain
\begin{equation}\label{xd12sfgh600}
\frac{d}{dt} g^t_* \omega = g^t_* (d i_v \omega + i_v d \omega) = d(g^t_* i_v \omega),
\end{equation}
where we used the fact that $d\omega = 0$ and that the exterior derivative commutes with the pullback ($g^t_* d = d g^t_*$).

Integrating both sides of \eqref{xd12sfgh600} with respect to $t$ from $-\infty$ to $0$, we get:
$$\int_{-\infty}^0 \frac{d}{dt} (g^t_* \omega) dt = \int_{-\infty}^0 d(g^t_* i_v \omega) dt.$$
Since $\lim_{t \to -\infty} g^t_* \omega = 0$ (see (\ref{dfh700})), the left-hand side equals $g^0_* \omega = \omega$. By pulling the exterior derivative out of the integral, we obtain:
$$\omega = d \left( \int_{-\infty}^0 g^t_* i_v \omega \, dt \right).$$
Thus, $\Omega = \int_{-\infty}^0 g^t_* i_v \omega \, dt$ is the required $(k-1)$-form.

All formulas involving improper integrals are justified by the following asymptotic relations:
$$g^t_* i_v \omega=O(e^{kt}),\quad d(g^t_* i_v \omega)=O(e^{kt})$$
holding as $t\to-\infty$.

An alternative proof is given in \cite{SCHwartz}.\end{proof}

\end{document}